\documentclass[11pt,a4]{article}
\usepackage{graphicx}
\usepackage{amsmath}
\usepackage[numbers]{natbib}
\usepackage{pslatex}
\usepackage{color}
\usepackage{SPE}
\usepackage{algorithm}
\usepackage{algpseudocode}
\usepackage{hyperref}
\usepackage{multirow}
\usepackage{bm}
\usepackage{colortbl}

\begin{document}
%Title
\title{Decoupled Block-Wise ILU(k) Preconditioner on GPU}

\author{Bo Yang, Hui Liu, He Zhong and Zhangxin Chen\\
University of Calgary}
\date{}
\maketitle

{\noindent \rule{\linewidth}{0.5mm}}

\bigskip

\thispagestyle{empty}

\section*{Abstract}
This research investigates the implementation mechanism of block-wise ILU(k) preconditioner on GPU. The block-wise ILU(k) algorithm requires both the level $k$ and the block size to be designed as variables. A decoupled ILU(k) algorithm consists of a symbolic phase and a factorization phase. In the symbolic phase, a ILU(k) nonzero pattern is established from the point-wise structure extracted from a block-wise matrix. In the factorization phase, the block-wise matrix with a variable block size is factorized into a block lower triangular matrix and a block upper triangular matrix. And a further diagonal factorization is required to perform on the block upper triangular matrix for adapting a parallel triangular solver on GPU. We also present the numerical experiments to study the preconditioner actions on different $k$ levels and block sizes.

\section*{\large Keywords}
ILU (k), block-wise matrix, parallel computing, GPU, preconditioner

\section{Introduction}

When running a reservoir simulator, over 70\% of time is spent on the solution of linear systems derived from the Newton methods \cite{CHM}. If large and highly heterogeneous geological models are applied, the linear systems are even harder to solve and require far more simulation time. The Krylov iterative linear solvers and preconditioners for the solution of linear systems have been studied for decades. The GMRES (Generalized Minimal Residual) method and BiCGSTAB (Biconjugate Gradient Stabilized) method are general algorithms for linear systems \cite{template,Y}. The ORTHOMIN (orthogonal minimum residual) method is originally designed by Vinsome for reservoir simulations \cite{orthomin}. Both the robustness and efficiency of these iterative solvers can be optimized by using preconditioners. Many advanced preconditioners have been developed, such as the domain decomposition preconditioners \cite{cai}, multi-stage preconditioners \cite{CHM, jcxu}, a constrained pressure residual (CPR) preconditioner \cite{CPR-cao, NCC} and a fast auxiliary space preconditioner (FASP) \cite{jcxu}. The ILU is one of the most commonly used preconditioners \cite{template, Y, CZ_dant, LS_gpil, KSLS}.

In reservoir simulation, block-wise matrices are frequently encountered. For example, a reservoir model consists of many grid blocks in a saturated black oil case. Each block has three unknowns (the pressure of the oil phase, the saturation of the water phase and the saturation of the gas phase). If the unknowns in each block are numbered consecutively, the discrete system of the mass equilibrium equations from the Newton method is a block-wise matrix where each block is a $3\times3$ sub matrix. For other reservoir models, such as a compositional model and a thermal model, more unknowns are in a grid block. Because the block size is determined by the number of unknowns, it varies and depends on a concrete model. Therefore, a preconditioner algorithm must be compatible to block matrices with a variable block size.

GPUs (Graphical Processing Units) are special electronic devices designed to manipulate and accelerate the creation of frame images displaying on a screen. Large amounts of pixel data can be computed on GPUs in parallel. Nowadays, the applications of GPUs have been more and more popular in scientific computing fields. Modern GPUs own high memory speed and floating point performance. For example, the memory speed of NVIDIA Tesla K40 can reach 288 GB/s. It provides up to 1.66 Tflops in double precision floating point calculations. In contrast, the memory speed of CPUs (Central Processing Units) is relatively low. For example, the theoretical peak bandwidth of a DDR4-3200 memory
is 51.2 GB/s on a 128-bit bus. The peak performance of an Intel Core i7-5960X Processor Extreme Edition CPU is only 385 Gflops in double precision. GPUs and CPUs are designed for different purposes. CPUs are utilized in solving serial and complicated tasks while GPUs are specialized in parallel computing. GPUs should be more propriety when an algorithm has enormous parallel components. Some solvers for linear systems have been developed on GPUs, such as the Krylov linear solvers \cite{NCC, KSLS, LS, CLY}. Haase et al. developed a parallel AMG (algebraic multigrid) solver using a GPU cluster \cite{HLDP}. Chen et al. designed classical AMG solvers on GPU \cite{uc-amg}. Bell et al. investigated fine-grained parallelism of AMG solvers on GPU \cite{amg-nv}. Bolz, Buatois, Goddeke, Bell, Wang, Brannick, Stone and their collaborators also studied GPU-based parallel algebraic multigrid solvers, and details can be found in \cite{bolz, bua, god, bran, nbell, lwang}. NVIDIA developed a hybrid matrix format HYB and sparse Krylov subspace solvers \cite{NCC}, which are used for general sparse matrices \cite{nv-spmv, nv-spmv2}. NVIDIA also implemented the scientific computing libraries FFT (fast Fourier transform) \cite{guide} and BLAS (basic linear algebra subprograms) \cite{nv-spmv, nv-spmv2, guide}.

We have studied and developed a set of block-wise ILU (k) preconditioned Krylov solvers \cite{block_iluk}. Because of the high complexity of block-wise ILU (k) implementation, we describe the detailed mechanism of the decoupled block-wise ILU (k) preconditioner in this technical report. The main ideas are stated below. 

We assume $A$ represents the coefficient matrix of a linear system. A preconditioner system $Mx = b$ is supposed to be solved at least once in each iteration of a Krylov solver. Based on ILU, the matrix $M$ is expressed as $LU$, where $L$ is a lower triangular matrix and $U$ is an upper triangular matrix. Then the preconditioner system can be easily decomposed into $Ly = b$ and $Ux = y$ which would be solved by a triangular solver. The ILU factorization derives from Gauss elimination by applying a nonzero pattern $P$ as a factorization filter. The ILU factorization is rather cheap to compute when $P$ is sparse. If $P$ has the same nonzero pattern as $A$, the factorization is called ILU(0). Because $A$ is sparse, sometimes the pattern cannot provide enough nonzero positions for $L$ and $U$, which reduces the rate of convergence. An improved approach is ILU(k) that allows additional fill-in positions on $P$ and more accurate factorizations can be obtained. The level $k$ is designed to control the extent of fill-in. A higher $k$ means more fill-in.

The ILU(k) algorithm consists of two steps. One is the creation of the fill-in nonzero pattern $P$. The other is the ILU factorization based on the pattern $P$. We call the first step a symbolic phase and the second step a factorization phase. Because the procedure of establishing $P$ has nothing to do with the concrete entry values of the original matrix $A$, a better way is to separate the symbolic phase from the ILU(k) algorithm and design it solely. After extracting the symbolic phase, the rest of the algorithm looks like ILU(0), which uses the fill-in nonzero pattern $P$ to perform an ILU factorization on the original matrix $A$. The ILU(0) algorithm requires that the nonzero pattern comes from $A$. In our algorithm, matrix $A$ is saved by data structure CSR (Compressed Sparse Row) which only contains the nonzero entries and their positions. If we use zero values to fill the positions of fill-in positions in matrix $A$, it will have the same non-zero pattern as $P$ without any essential value changes. Thus, the ILU(0) algorithm can be applied to a modified $A$ directly.

The linear systems we will solve are block-wise systems. A block ILU(k) preconditioner is a better choice than a point-wise ILU(k) preconditioner in general especially when the condition number of $A$ is large. The block-wise preconditioner is more stable even if the matrix is ill-conditioned. Because the block size depends on the concrete mathematical models from reservoirs, the block-wise ILU(k) should be compatible to any block size. The main idea of applying ILU(k) to a block-wise matrix is to extract a point-wise matrix structure from the block-wise matrix $A$, where each nonzero block of $A$ is looked at as a nonzero entry. The fill-in nonzero pattern $P$ is created from the point-wise structure by ILU(k). Then the pattern $P$ is used to factorize $A$ into block $L$ and block $U$ by operating on each block as a whole. Obviously, the whole ILU(k) algorithm for the block-wise matrix $A$ mixes up point-wise matrix operations with block-wise matrix operations. As stated above, we already decoupled the symbolic phase and the factorization phase of the point-wise ILU(k). The same idea can be applied in the block-wise ILU(k) factorization. Thus, we can focus on the creation of the fill-in nonzero pattern $P$ according to the level $k$ in the symbolic phase and then on producing $L$ and $U$ for a block-wise matrix in the factorization phase.

A block preconditioner system can be expressed as $LUx = b$, where $L$ and $U$ are both block matrices generated by ILU(k). This equation can be decomposed to two block triangular equations $Ly = b$ and $Ux = y$. However, the parallelization of the block ILU(k) is difficult to implement on a NVIDIA GPU platform. Because the calculation data is stored on a NVIDIA GPU by the global memory \cite{guide}, the data access pattern on GPUs must be coalesced to achieve optimal efficiency. It is hard to compute a coalesced access size if the block size of a matrix is variable. In our research, we convert block triangular systems to point-wise triangular systems to solve. Thus a point-wise parallel triangular solver on GPU can be used. Nevertheless, a new problem emerges; the point-wise parallel triangular solver requires that the matrix must be a point-wise triangular matrix. The block matrix $L$ meets the requirement because its diagonal blocks are all identity blocks. The block matrix $U$ does not satisfy the criterion because it cannot guarantee all the diagonal blocks are diagonal matrices. To overcome this problem, a further decomposition for the block matrix $U$ is required. Based on $U = D(D^{-1}U)$, a block diagonal matrix is extracted from $U$ which forms an additional system $z = D^{-1}y$ to be solved by a matrix-vector multiplication. The remaining $D^{-1}U$ is a block upper triangular matrix with pure identity blocks in its diagonal, which can be solved by a point-wise triangular solver on GPU. In our numerical experiments, the decoupled block ILU(k) algorithm runs smoothly with different configurations of $k$ and a block size. It shows favorable speedups against a traditional CPU. When the parameter $k$ or block size increases, the speedup effect may decrease. The reasoning analysis is also detailed in the numerical experiment section.

The layout of this paper is as follows:
In \S \ref{sec-blockwiseILUk}, we describe the detailed implementation mechanism of the decoupled block-wise ILU (k) preconditioner.
In \S \ref{sec-experiments}, numerical experiments are performed to test the speedup performance on GPU.
In \S \ref{sec-conclusion}, some conclusions are given.

%-------------------------------------------------------------------------
\section{Decoupled Block-wise ILU(k) Preconditioner}
\label{sec-blockwiseILUk}
\subsection{Data structure }
\label{subsec-Datastructure}

\begin{figure}[tbh]
    \centering
    \includegraphics[width=0.60\linewidth]{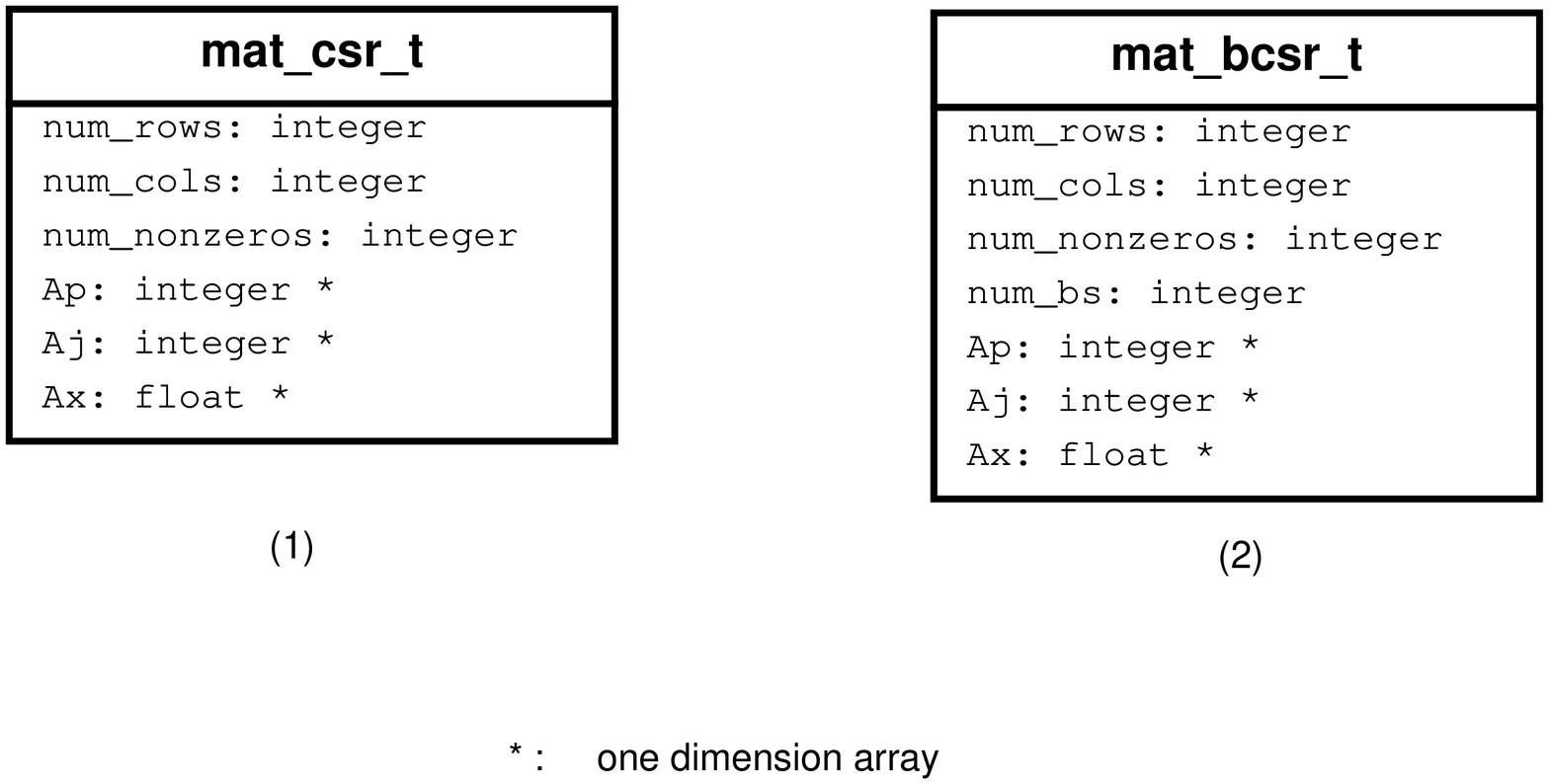}
    \caption{Data structure for matrices.}
    \label{fig-mat_csr_t}
\end{figure}

\begin{figure}[tbh]
    \centering
    \includegraphics[width=0.85\linewidth]{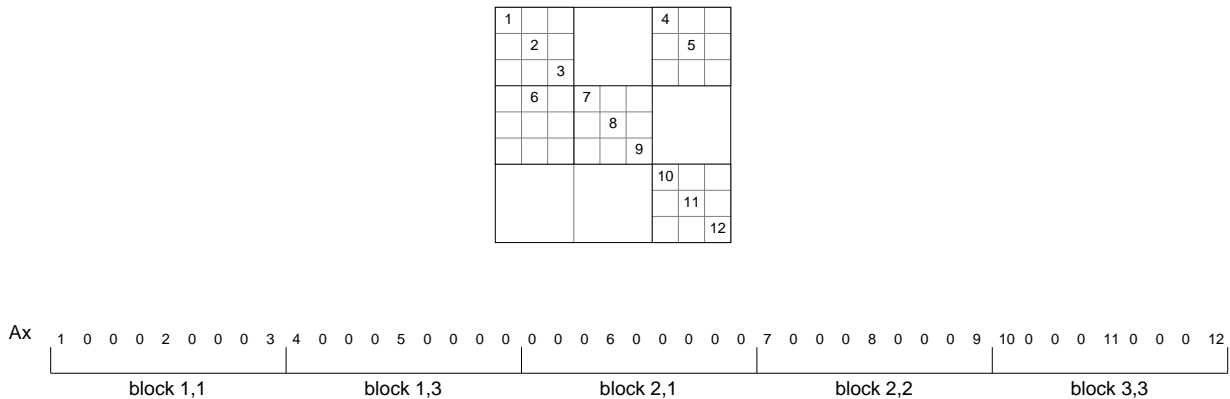}
    \caption{Data storage approach for $mat\_bcsr\_t$.}
    \label{fig-Ax}
\end{figure}

In our study, all matrices are stored in the CSR (Compressed Sparse Row) format. Data structures $mat\_csr\_t$ and $mat\_bcsr\_t$ are designed
for a point-wise matrix and a block-wise matrix, respectively; see Figure~\ref{fig-mat_csr_t}. The CSR format only keeps the nonzero elements by rows.
All nonzero elements are stored in the property $Ax$ row by row. The property $Aj$ is used for recording the column indices corresponding to each nonzero element in $Ax$.
The property $Ap$ records the start index of each row in array $Aj$ or $Ax$. For a block-wise matrix, all blocks can be classified into two categories that
are zero blocks and nonzero blocks. All the values of the elements in a zero block are zero, but a nonzero block contains at least one nonzero element.
If each zero block is looked at as a zero element and each nonzero block is looked at as a nonzero element, the whole block-wise matrix can be abstracted
as a point-wise matrix structure. In data structure $mat\_bcsr\_t$, the property $num\_bs$ represents the block size, and the properties $num\_rows$, $num\_cols$, $num\_nonzeros$,
$Ap$ and $Aj$ are all based on the concept of a nonzero block element. All the nonzero block elements are stored in $Ax$ row by row. The size of array $Ax$
is the number of nonzero blocks multiplied by the block size square. All elements in a non-zero block are saved, no matter if it is a zero element.
To keep compatible with BLAS (Basic Linear Algebra Subprograms) and LAPACK (a mathematical library of linear algebra routines for dense systems solution
and eigenvalue calculations), the elements of each block are stored in $Ax$ column by column. Figure~\ref{fig-Ax} gives a schematic of storing data for a block-wise matrix.
Both $mat\_csr\_t$ and $mat\_bcsr\_t$ use $Ap$ and $Aj$ to describe the structure of a matrix. When establishing an abstract point-wise matrix for a block-wise matrix,
$Ap$ and $Aj$ of $mat\_bcsr\_t$ can be used directly as the structure of the abstract point-wise matrix. The only difference is that each element in a block-wise matrix
is a block element, but each element in an abstract point-wise matrix is a point element. If ignoring the concrete values of elements and only considering whether
they equal zero or not, all the nonzero positions of a matrix form a nonzero pattern. An element is a nonzero element if its row and column indices are
involved in $Ap$ and $Aj$; otherwise, it is a zero element. Therefore, $Ap$ and $Aj$ also represent the nonzero pattern of a matrix.

\begin{figure}[tbh]
    \centering
    \includegraphics[width=0.20\linewidth]{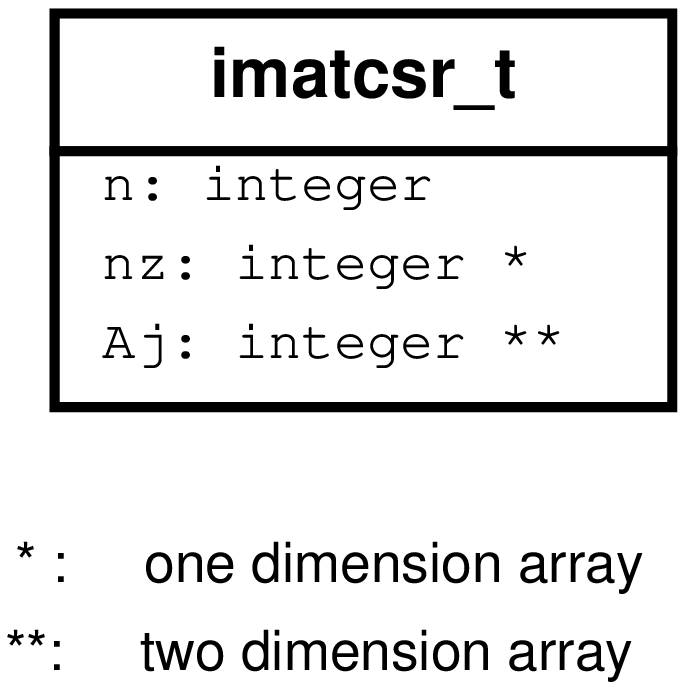}
    \caption{Data structure for nonzero pattern.}
    \label{fig-imatcsr_t}
\end{figure}

Figure~\ref{fig-imatcsr_t} shows data structure $imatcsr\_t$ which is specially designed for storing a nonzero pattern and will be applied in the symbolic
phase of ILU(k). $n$ represents the dimension of a matrix. $nz$ is an array for the length of each row. The two-dimensional array $Aj$ is used for storing column indices.
As stated above $mat\_csr\_t$ stores column indices in a one-dimensional array. If an element is deleted from the matrix, the entire elements after the current
element in the memory space must be moved forward, which costs a large amount of computing resources, especially when the deletion operation happens frequently. However, because the column indices are stored in the two-dimensional array $Aj$ of $imatcsr\_t$ row by row and each row contains sparse elements, a very low cost of moving data is required for the deleting operation. Therefore, $imatcsr\_t$ is more suitable for manipulating a nonzero pattern.

\subsection{Decoupled ILU(k) algorithm}
\label{subsec-DecoupledILUkalgorithm}

\begin{algorithm}[htb]
\caption{Point-wise ILU(0) factorization} \label{alg_1}
\begin{algorithmic}[1]
 \For {i = 2: n}
  \For {p = 1: i - 1 \& $(i, p) \in P$}
   \State $a_{ip} = a_{ip} / a_{pp}$
     \For {j = p + 1: n \& $(i, j) \in P$}
     \State $a_{ij} = a_{ij} - a_{ip}a_{pj}$
     \EndFor
  \EndFor
\EndFor
\end{algorithmic}
\end{algorithm}

\begin{figure}[tbh]
    \centering
    \includegraphics[width=0.50\linewidth]{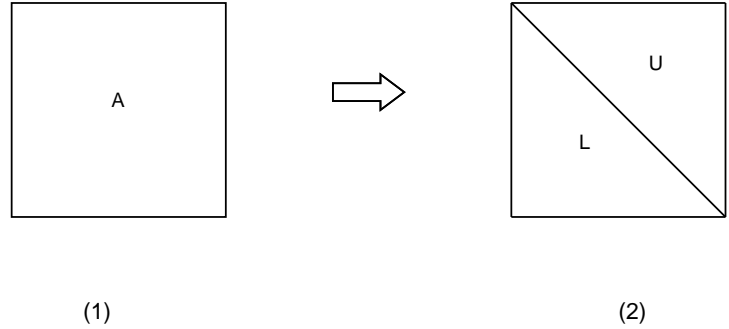}
    \caption{ILU factorization.}
    \label{fig-ILU}
\end{figure}

Algorithm~\ref{alg_1} shows the point-wise ILU(0) algorithm. $P$ represents the nonzero pattern of matrix $A$. If the code $(i,k) \in P$  on lines 2 and 4 is removed,
this algorithm is Gauss elimination. As we know the Gauss elimination is costly on computing resources. In the ILU(0) algorithm, $P$ serves as a filter
where only the elements lying in the nonzero positions of $P$ can be calculated by lines 3 and 5. Due to the sparse character of $P$, the complexity
of this algorithm is very low, which can be performed at each iteration of a Krylov subspace solver. Matrix $A$ is stored in $mat\_csr\_t$ format.
In this algorithm, the resulting matrix also uses the same memory as $A$, which means that the original elements of $A$ will be replaced by the resulting elements
calculated by lines 3 and 5. The lower triangular matrix $L$ from factorization will be stored in $A$'s low triangular part except the diagonal
which are all unit values. $L$ is generated by line 3. The upper triangular matrix $U$ from factorization will be stored in $A$'s upper triangular
part including the diagonal. $U$ is generated by line 5. Figure~\ref{fig-ILU} shows a schematic of ILU(0) factorization. Because $P$ comes from the nonzero pattern
of $A$, Figure~\ref{fig-ILU}-(2), which represents the resulting matrix, has the same nonzero pattern as Figure~\ref{fig-ILU}-(1).

Because $A$ is a sparse matrix, the nonzero pattern $P$ of $A$ is also sparse. The accuracy of ILU(0) may be insufficient to yield an adequate rate
of convergence. More accurate incomplete LU factorizations are often more efficient as well as more reliable. These more accurate factorizations differ from ILU(0) by allowing some fill-in. A level of $k$ is defined to control the degree of fill-in \cite{Y}. A larger $k$ allows more fill-ins
in the nonzero pattern in addition to the original nonzero pattern created by $A$. $k$ is a nonnegative number; ILU(k) is ILU(0) when $k$ is zero.

\begin{equation}
\label{equ-1}
L_{ij} = \left\{
 \begin{aligned}
 & 0, & (i, j) \in P \\
 & \infty, & (i, j) \notin P.
 \end{aligned}
  \right.
\end{equation}

\begin{equation}
\label{equ-2}
L_{ij} = min \{L_{ij}, L_{ip} + L_{pj} + 1\}.
\end{equation}

ILU(k) requires that a level is defined for each entry of matrix $A$. The initial level of each entry $A_{ij}$ is defined by formula~\eqref{equ-1},
where $P$ represents the nonzero pattern of matrix $A$. The nonzero elements of $A$ have the level zero; otherwise, it has a level of infinity.
Formula~\eqref{equ-2} gives the level updated algorithm \cite{Y}. Apparently, formula~\eqref{equ-2} has no influence on the level of a nonzero element, because the level of a
nonzero element always remains zero. This means that only zero element's level can be updated from infinity to a limited positive value. If the new level value
is still greater than level $k$, this entry will be kept as zero and the position of this entry is not involved in the next factorization.
In other words, only the positions that have a new level value less than level $k$ can continue to a next calculation. If this position continues
to the end, it will be a fill-in position included in the extended nonzero pattern. Because the nonzero pattern generated by ILU(k) has more nonzero positions than ILU(0),
the accuracy of solution for a preconditioner system improves. Algorithm~\ref{alg_2} gives the whole procedure of ILU(k).

\begin{algorithm}[htb]
\caption{Point-wise ILU(k) factorization} \label{alg_2}
\begin{algorithmic}[1]
\State For all nonzero entries in nonzero pattern $P$ define $L_{ij} = 0$
 \For {$i = 2: n$}
  \For {$p = 1: i - 1 $ \& $L_{ip} \le k$}
   \State $a_{ip} = a_{ip} / a_{pp}$
     \For {$j = p + 1: n $}
     \State $a_{ij} = a_{ij} - a_{ip}a_{pj}$
     \State $L_{ij} = min \{L_{ij}, L_{ip} + L_{pj} + 1\}$
     \EndFor
  \EndFor
  \If {$L_{ij} > k$}
     \State $a_{ij} = 0$
  \EndIf
\EndFor
\end{algorithmic}
\end{algorithm}

Apparently, the calculations of $a_{ip}$, $a_{ij}$ and $L_{ij}$ are mixed together in Algorithm~\ref{alg_2}, which complicates the algorithm. There is no direct
relationship between $a_{ip}$, $a_{ij}$ and $L_{ij}$. $L_{ij}$ merely acts as a filter to control whether $a_{ip}$ or $a_{ij}$ is qualified to be calculated
in the following steps. Because different functions are compounded into one subroutine to implement, the design of Algorithm~\ref{alg_2} decreases the maintainability
of the codes. From another perspective, line 6 is not necessary to be calculated when $L_{ij}$ is greater than $k$ because $a_{ij}$ is set to zero on line 11.
In other words, if an element $a_{ij}$ cannot stay to the end to be a fill-in element, the calculation for $a_{ij}$ on lines 6 is unnecessary. Thus line 6 wastes
some computing resources in some situations. Note that our purpose is to implement ILU(k) for both variable level $k$ and block size. Level $k$ relates to
the creation of a nonzero pattern and the block size relates to the ILU factorization. It is more clear and easier to process if Algorithm~\ref{alg_2}
is decoupled instead of implementing everything in one subroutine. For these reasons, if the functions of the zero pattern creation and the ILU factorization can be
separated, all the problems stated above can be resolved. Therefore, for both the aspects of designing a favorable maintainable program and a performance optimization
algorithm, a decoupled ILU(k) implementation is necessary. The decoupled ILU(k) contains two steps. The first step is called the symbol phase focusing on establishing
a fill-in nonzero pattern, which has no relation to the processing of any concrete data values. The second phase is responsible for the ILU factorization which uses the nonzero
pattern established in the symbolic phase to factorize the original matrix into a lower triangular part and an upper triangular part.

\begin{algorithm}[htb]
\caption{Symbolic phase} \label{alg_3}
\begin{algorithmic}[1]
\State For all nonzero entries in nonzero pattern $P$ define $L_{ij} = 0$
\State Define $P'$ as $n \times n$ nonzero pattern
\State Initiate $P'$ by full filling with entries
 \For {$i = 2: n$}
  \For {$p = 1: i - 1 $ \& $L_{ip} \le k$}
     \For {$j = p + 1: n $}
     \State $L_{ij} = min \{L_{ij}, L_{ip} + L_{pj} + 1\}$
     \EndFor
  \EndFor
  \If {$L_{ij} > k$}
     \State remove entry $ij$ from $P'$
  \EndIf
\EndFor
\end{algorithmic}
\end{algorithm}

The symbolic phase can be presented in Algorithm~\ref{alg_3}, where $P'$ is designed to record the fill-in nonzero pattern. The data structure $tm\_ilu\_t$ is applied to $P'$,
which can easily and efficiently remove entries. The logic for computing $L_{ij}$ is the same as in Algorithm~\ref{alg_2}. Algorithm~\ref{alg_3} has been separated totally
from the processing of concrete data values. The factorization phase is responsible for creating $L$ and $U$ according to the nonzero pattern $P'$. Because $tm\_ilu\_t$ only stores the pattern without actual data values, in order to apply the ILU factorization, it is necessary to create a matrix $A'$ by $mat\_csr\_t$ from the nonzero pattern $P'$, which contains both structure and data values. All entries
are copied from $A$ to $A'$ since all the entry positions of $A$ are remained in $P'$. The fill-in positions of $A'$ are all filled with a value zero. Thus $A'$ has the same
structure as $P'$ but has no essential data value difference from $A$. The ILU(0) algorithm can be directly applied to $A'$ so it is called the factorization phase; see
Algorithm~\ref{alg_1}. The process of the factorization phase is also illustrated in Figure~\ref{fig-ILU}. The only difference from ILU(0) for $A$ is that the number of nonzero
elements in Figure~\ref{fig-ILU}-(2) is greater than the number of nonzero elements in $A$, because extra fill-in elements are imported by ILU(k).

\subsection{Decoupled ILU(k) for block-wise matrix}
\label{subsec-DecoupledILUkforblockwisematrix}

\begin{equation}
\label{equ-3}
A = \left(
    \begin{array}{cccc}
    A_{11} & A_{12} & \cdots  & A_{1n} \\
    A_{21} & A_{22} & \cdots  & A_{2n} \\
    \cdots & \cdots & \cdots  & \cdots \\
    A_{i1} & A_{i2} & \cdots  & A_{in} \\
    \cdots & \cdots & \cdots  & \cdots \\
    A_{n1} & A_{n2} & \cdots  & A_{nn}
    \end{array}
    \right)
\end{equation}
As mentioned above, for reservoir simulations, each grid block has several unknowns, such as pressure, temperature and saturation. If all unknowns in each block are numbered consecutively,
matrix $A$ from the Newton methods is a block-wise matrix and has the structure as formula (3), where $A$ is a nonsingular block matrix, each block $A_{ij}$  $(1 \le i,j \le n)$ is an $m \times m$ submatrix, and $m$ represents the block size. Obviously, the block matrix is a regular point-wise matrix when $m$ is one. The block matrix $A$ is also equivalent to a $nm \times nm$ point-wise matrix from a point-wise view. Thus the Krylov subspace solvers with the point-wise ILU(k)
preconditioner can be employed to solve the block-wise system. However, for a block-wise matrix, a block ILU(k) preconditioner can be a better choice than a point-wise ILU(k)
preconditioner when the condition number of $A$ is large. The point- and block-wise preconditioners have been studied and implemented
on CPUs by Saad et al. \cite{Y}. Our implementation mechanism of a decoupled block ILU(k) preconditioner on GPU is described below.

\begin{figure}[tbh]
    \centering
    \includegraphics[width=0.80\linewidth]{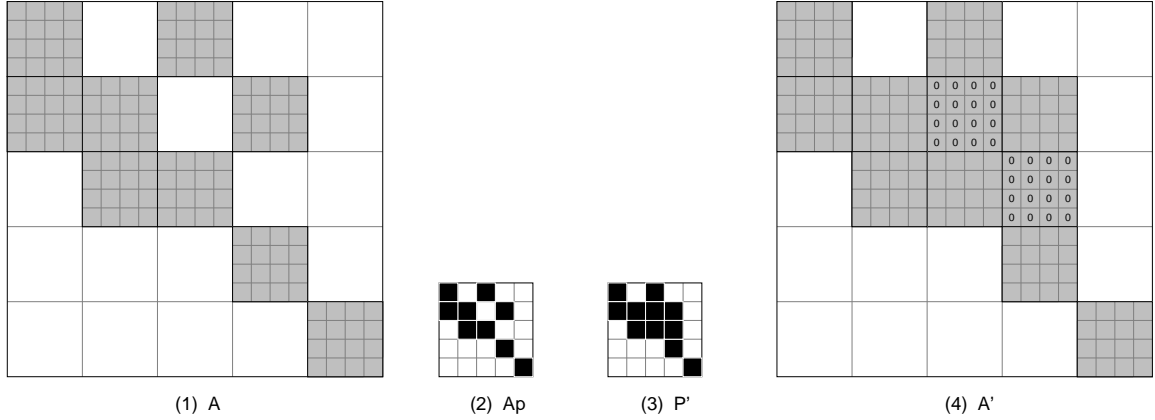}
    \caption{Block-wise ILU(k) factorization.}
    \label{fig-symbol_pattern}
\end{figure}

First, we need to abstract matrix $A$ into a point-wise matrix $Ap$. For matrix $A$, if block $A_{ij}$ has at least one nonzero element, block $A_{ij}$ is called a
nonzero block. All nonzero blocks are looked at as a nonzero element in $Ap$. The symbolic phase algorithm is applied to $Ap$ to create a fill-in nonzero pattern. Thus
the entries in $Ap$ do not need to have any concrete data values and $Ap$ can be stored by format $imatcsr\_t$. Figures~\ref{fig-symbol_pattern}-(1) and~\ref{fig-symbol_pattern}-(2)
give a schematic of such a process. Figure~\ref{fig-symbol_pattern}-(1) represents the block-wise matrix $A$ where the gray blocks denote the nonzero blocks.
Figure~\ref{fig-symbol_pattern}-(2) represents the point-wise matrix $Ap$ where each black point is a nonzero position extracted from $A$. Second, the fill-in nonzero
pattern $P'$ is established in Algorithm~\ref{alg_3} on $Ap$, see Figure~\ref{fig-symbol_pattern}-(3). The next step is to create the block-wise matrix $A'$ used for
the factorization phase according to the pattern $P'$. All nonzero blocks of $A'$ are copied from $A$. All the fill-in positions are filled by zero blocks. The process is shown in
Figures~\ref{fig-symbol_pattern}-(3) and~\ref{fig-symbol_pattern}-(4). Apparently, Figure~\ref{fig-symbol_pattern}-(3) has more fill-in positions than Figure~\ref{fig-symbol_pattern}-(2)
and Figure~\ref{fig-symbol_pattern}-(4) has more fill-in zero blocks than Figure~\ref{fig-symbol_pattern}-(1). Again, the block matrix $A'$ and the block matrix $A$ have no difference
in data values. However, from the matrix structure point of view, $A'$ conforms to the nonzero pattern $P'$ and can be used for a block-wise ILU(0) factorization directly.
The algorithm of ILU(0) for a block-wise matrix is a little different from the algorithm of ILU(0) for a point-wise matrix, which is shown in Algorithm~\ref{alg_4}.
All operations on matrix $A'$ are matrix operations on blocks, such as the inverse of the diagonal block $A_{pp}$  and the matrix-matrix multiplication $A_{ip}A_{pj}$.
The block-wise ILU(k) factorization result is shown in Figure~\ref{fig_ILU0}-(1). The integrated decoupled block ILU(k) algorithm is given in Algorithm~\ref{alg_5}.

\begin{algorithm}[htb]
\caption{Block-wise ILU(0) factorization algorithm} \label{alg_4}
\begin{algorithmic}[1]
 \For {i = 2: n}
  \For {p = 1: i - 1 \& $(i, p) \in P'$}
   \State $A_{ip} = A_{ip}A_{pp}^{-1}$
     \For {j = p + 1: n \& $(i, j) \in P'$}
     \State $A_{ij} = A_{ij} - A_{ip}A_{pj}$
     \EndFor
  \EndFor
\EndFor
\end{algorithmic}
\end{algorithm}

\begin{algorithm}[htb]
\caption{Decoupled block ILU(k) algorithm} \label{alg_5}
\begin{algorithmic}[1]
\State Block-wise matrix $A$ is abstracted as point-wise matrix $Ap$
\State Fill-in nonzero pattern $P'$ is established through symbolic phase on
      $Ap$ by Algorithm~\ref{alg_3}
\State Backfill data of $A$ to block-wise matrix $A'$
\State Fill zero blocks in $A'$ according to the fill-in positions from $P'$
\State Block-wise ILU(0) factorization on $A'$ by Algorithm~\ref{alg_4}
\end{algorithmic}
\end{algorithm}

\subsection{Parallel triangular solver on GPU}
\label{subsec-ParalleltriangularsolveronGPU}

After assembling the block triangular matrices $L$ and $U$, the precoditioner system $LU = b$ are required to solve. A parallel triangular solver can be used to
solve triangular systems on GPU. The principle of a parallel triangular solver is the level schedule method \cite{Y, LS_gpil}. It has been developed by Chen et al. \cite{ZHB-ptsg}.
The idea is to group unknowns $x_i$ into different levels so that all unknowns within the same level can be computed simultaneously \cite{Y, LS_gpil}. Because an upper triangular system can be easily transform to a lower triangular system, only a lower triangular system is analyzed below. The level of $x_i$ can be defined as in formula~\eqref{equ-4}. The level schedule method
is described in Algorithm~\ref{alg_6}. For GPU computing, each level in Algorithm~\ref{alg_6} can be parallelized:
\begin{equation}
\label{equ-4}
 l(i) = 1 + \max_j {l(j)} \quad \text{ for  all } j \text{ such
that } \ L_{ij} \neq 0, i = 1, 2, \ldots, n,
\end{equation}
where
\begin{itemize}\itemsep1pt \parskip0pt \parsep0pt
  \item $L_{ij}$ : the $(i, j)$th entry of $L$
  \item $l(i)$ : the $i$th level, zero initially
  \item $n$ : the number of rows
\end{itemize}

\begin{algorithm}[htb]
\caption{Level schedule method for solving a lower triangular system}
\label{alg_6}
\begin{algorithmic}[1]
\State \textbf{input}
\State n			         \Comment the number of levels
\For {i = 1 : n}
	\For {each row in current level}
		\State solve the current row
	\EndFor
\EndFor
\end{algorithmic}
\end{algorithm}

\begin{figure}[tbh]
    \centering
    \includegraphics[width=0.80\linewidth]{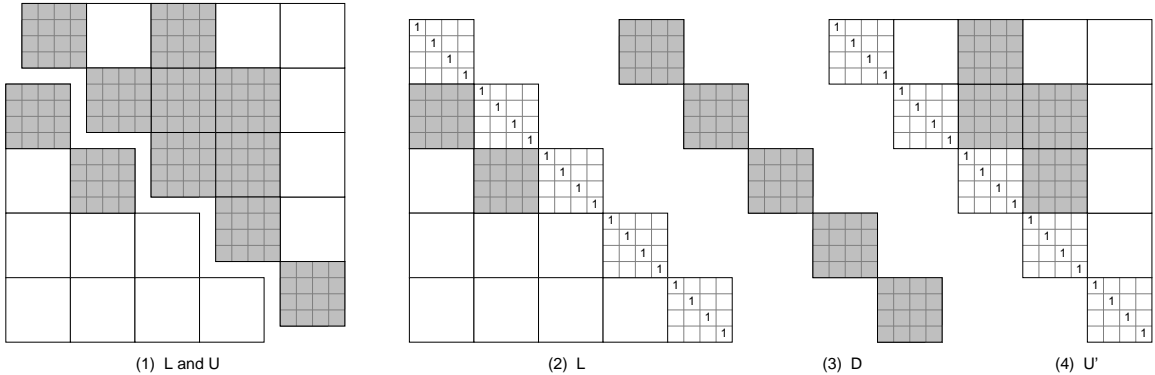}
    \caption{Block upper triangular matrix factorization.}
    \label{fig_ILU0}
\end{figure}

Though the level schedule method is designed for a point-wise matrix, it can be directly used for the block lower triangular matrix $L$. Because the diagonal blocks in $L$
are all identity matrix, $L$ can be looked at as a point-wise lower triangular matrix straightly; see Figure~\ref{fig_ILU0}-(2). However, the block upper triangular matrix $U$
is not an upper triangular matrix from a point-wise view because its diagonal blocks are not always point-wise diagonal matrices. The diagonal of $U$ presents a
kind of irregular zigzag shape. To resolve this issue, $U$ can be factorized into a multiplication of a block diagonal matrix $D$ and a unit diagonal
matrix $U'$; see equation~\eqref{equ_5}. A schematic is given in Figures~\ref{fig_ILU0}-(3) and~\ref{fig_ILU0}-(4).

\begin{equation}
\label{equ_5}
U=D(D^{-1} U)=DU'
\end{equation}

\begin{equation}
\label{equ_6}
Ly = b,z = D^{-1} y,U' x = z
\end{equation}

The linear systems $Ly = b$ and $U'x = z$ can be solved by the point-wise level schedule method, and the system  $z = D^{-1}y$ is a direct matrix-vector multiplication operation.
The sparse matrix vector multiplication on GPU has been developed by Chen et al. \cite{LYCHS_spmv}. By far, the parallel solver for block triangular linear systems
on GPU is developed, which is appropriate for any block size. The algorithm is summarized in Algorithm~\ref{alg_7}.

\begin{algorithm}[htb]
\caption{Parallel solver for $LDU'x = b$}
\label{alg_7}
\begin{algorithmic}[1]
\State Solve $Ly = b$ by Algorithm~\ref{alg_6}
\State Solve $z = D^{-1}$y by sparse matrix-vector multiplication
\State Solve $U'x = z$ by Algorithm~\ref{alg_6}
\end{algorithmic}
\end{algorithm}

\section{Numerical Results}
\label{sec-experiments}

In this section, a series of numerical experiments are employed to test the performance of the block-wise ILU(k) preconditioner on GPU against a CPU.
Our workstation is composed of an Intel(R) Xeon(R) CPU E5-2680 0 @ 2.70GHz CPU and a NVIDIA Tesla K20Xm GPU which has a 249.6 GB/s memory bandwidth and
2688 CUDA cores. The workstation operating system is the Red Hat Enterprise Linux Server release 6.6 (Santiago). The development environment is CUDA 6.5
and GCC 4.4. All floating data values are in double precision. The linear solver is GMRES with restarts 20. We also implement the same block ILU(k)
algorithm for a CPU, which plays as a contrastive algorithm. The CPU codes are all compiled with -O3 option and only one thread is employed.
An algorithm speedup on a GPU against a CPU is calculated by tc/tg where tc is the CPU time and tg is the GPU time.

The numerical experiments contain three parts which are performed for different matrices. Each experiment result is listed in an individual
table in Tables 1-3. Each table consists of several sections. A group of runs with a fixed
block size are organized in the same section. The effects of different $k$ levels can be judged within each section. Section one in each table represents
the point-wise ILU(k) where the block size equals one. The other three sections show block-wise ILU(k) effects for different block sizes. From a practical point of view,
the block sizes are determined by the underlying mathematical models. In the experiments, we set the block size to various values to compare the performance of the block-wise ILU(k)
preconditioner under different situations. There are many factors that affect the speed of solving a linear system, such as the number of fill-in elements,
the degree of parallelism, the ratio of nonzero elements to rows, and the pattern of nonzero elements. In this study, we only focus on the extent
of fill-in controlled by level $k$ and the degree of parallelism affected by the block size.

\textbf{Experiment 1}: The matrix used in this experiment is a 3D Poisson equation. Its dimension is 1,728,000 and has 12,009,600 nonzero elements.
The average number of nonzero elements on each row is seven. The results are displayed in Table~\ref{table-1}.

\begin{table}%[htb]
\centering
\begin{tabular}{cccccc} \hline

\bfseries Block size & \bfseries k level & \bfseries GPU time (second) & \bfseries CPU time (second) & \bfseries Iteration  & \bfseries Speedup  \\ \hline

1&	0&	3.8026&	    26.1277&	240&	6.9\\
1&	1&	3.2824&	    14.0474&	120&	4.3\\
1&	2&	7.3804&	    22.5474&	191&	3.1\\
1&	3&	5.4581&	    12.9243&	100&	2.4\\ \hline

2&	0&	3.2436&	    23.7607&	200&	7.3\\
2&	1&	2.7190&	    13.6358&	80&  	5.0\\
2&	2&	7.9980&	    22.9449&	157&	2.9\\
2&	3&	6.3791&	    11.6595&	70&	    1.8\\ \hline

4&	0&	3.4379&	    27.6770&	180&	8.1\\
4&	1&	3.9144&	    13.1689&	67&	    3.4\\
4&	2&	11.2787&	27.2849&	135&	2.4\\
4&	3&	11.1482&	16.5016&	65&	    1.5\\ \hline

8&	0&	4.6441&	    35.3785&	174&	7.6\\
8&	1&	6.7869&	    20.2972&	64&	    3.0\\
8&	2&	17.1339&	36.1713&	115&	2.1\\
8&	3&	15.3521&	24.4867&	60&	    1.6\\ \hline

\end{tabular}
\caption{Performance for the matrix 3D Poisson equation}
\label{table-1}
\end{table}

In Table~\ref{table-1}, the block-wise ILU(k) algorithm running time on both GPU and CPU is listed in columns GPU time and CPU time, respectively.
The speedup shows the computing speed ratio of the GPU time to the CPU time. All the speedup values are greater than one, which means that the algorithm speed on
a GPU is faster than on a CPU. In some runs, the speedup can be over five which is a favorable acceleration. The purpose of using the capability
of GPU parallel computing is realized. The speedup values have some fluctuations. The maximal value is 8.1 and the minimal value is 1.5.
Obviously, the speedup has a general decreasing tendency along with a higher $k$ level or a higher block size. The influence of
the $k$ level and the block size is analyzed according to the experiment data below.

Level $k$ of ILU(k) is used for adjusting the performance of convergence. More iterations reflect a lower convergence performance. We can see
the effects of different $k$ levels in each section. For example, the iteration number is 240 when the $k$ level is zero in section one. The iteration number reduces
to 100 when the $k$ level is 3. The general tendency of iterations is going down because a higher $k$ imports more fill-in positions, which leads to
a more accurate solution of solving a preconditioner system. However, fill-in positions also lead to two negative effects. One is the extra calculations
because of more elements imported. The other is a decrease in parallelism because of more relevance between rows caused by these fill-in elements.
This is reflected by the GPU time that shows a general increasing tendency with $k$ increasing in each section.

The influence of performance by the block size can be seen from the comparison between various sections. In general, a block-wise ILU(k) is a better preconditioner
than a point-wise ILU(k) for a block-wise matrix as mentioned above. A block-wise ILU(k) can run more steadily than a point-wise ILU(k) especially when a matrix
is ill-conditioned and the matrix conditioner number is large. However, a larger block size imports more fill-in zero elements in the blocks. Then the performance
must be affected by them. From Table 1, the GPU time ranges from 3.8026 to 7.3804 in section one where the block size is one. It increases to a new range which is
from 4.6441 to 17.1339 in section four where the block size is eight. Because a larger block size imports more zero elements in a block, extra computations are introduced.
In addition, more elements also reduce the degree of parallelism and lead to a slower GPU speed. Therefore, both reasons give a negative effect on GPU running and
prolong the GPU time.

\textbf{Experiment 2}: The matrix used in this experiment is matrix $parabolic\_fem$ which is downloaded from the University of Florida sparse matrix collection \cite{mmarket}.
Its dimension is 525,825 and has 2,100,225 nonzero elements. The average number of nonzero elements on each row is four. The results are listed in Table~\ref{table-2}.

\begin{table}%[htb]
\centering
\begin{tabular}{cccccc} \hline

\bfseries Block size & \bfseries k level & \bfseries GPU time (second) & \bfseries CPU time (second) & \bfseries Iteration  & \bfseries Speedup  \\ \hline

1&	0&	0.0655&	0.6186&	21&	9.4\\
1&	1&	0.0653&	0.6108&	21&	9.4\\
1&	2&	0.0655&	0.5935&	21&	9.1\\
1&	3&	0.0652&	0.5800&	21&	8.9\\ \hline

5&	0&	0.1222&	1.0407&	21&	8.5\\
5&	1&	0.1223&	1.0366&	21&	8.5\\
5&	2&	0.1224&	1.0387&	21&	8.5\\
5&	3&	0.1222&	1.0331&	21&	8.5\\ \hline

25&	0&	0.6849&	3.7078&	21&	5.4\\
25&	1&	0.6845&	3.7021&	21&	5.4\\
25&	2&	0.6861&	3.7376&	21&	5.4\\
25&	3&	0.6858&	3.7094&	21&	5.4\\ \hline

\end{tabular}
\caption{Performance for matrix $parabolic\_fem$}
\label{table-2}
\end{table}

All the speedup values range from 5 to 9 in Table~\ref{table-2}, which displays that the block-wise ILU(k) has higher performance on a GPU than a CPU.
In section one, the speedup is around nine. The speedup decreases slightly when $k$ increases. Because this matrix is very sparse and the ratio of nonzero
elements to rows is only 3.99, the fill-in positions increase slowly with a higher $k$. Therefore, the changing effects of speedup and iterations are not obvious.
Because the number of rows is 525,825 in the $parabolic\_fem$ matrix, the block sizes are designed as 1, 5 and 25 in this experiment. When the block size changes
from 1 to 5 and to 25, the speedup reduces from 9 to 8 and then to 5, respectively. At the same time, the GPU time also increases from around 0.065 to 0.685. These
effects show that the degree of parallelism is affected by the block size parameter.

\textbf{Experiment 3}: The matrix used in this experiment is matrix $atmosmodd$ which is downloaded from the University of Florida sparse matrix collection \cite{mmarket}.
Its dimension is 1,270,432 and has 8,814,880 nonzero elements. The average number of nonzero elements on each row is seven. The results are listed in Table~\ref{table-3}.

\begin{table}%[htb]
\centering
\begin{tabular}{cccccc} \hline

\bfseries Block size & \bfseries k level & \bfseries GPU time (second) & \bfseries CPU time (second) & \bfseries Iteration  & \bfseries Speedup  \\ \hline

1&	0&	1.5976&	     8.7331&	116&	 5.5\\
1&	1&	1.7882&	     6.1525&	77&	     3.4\\
1&	2&	3.5091&	     8.8263&	108&	 2.5\\
1&	3&	3.5308&	     7.0648&	78&	     2.0\\ \hline

2&	0&	1.3456&	     7.1068&	88&	     5.3\\
2&	1&	1.6299&	     6.0481&	58&	     3.7\\
2&	2&	3.0623&	     8.2178&	80&	     2.7\\
2&	3&	3.8071&	     6.8569&	58&	     1.8\\ \hline

4&	0&	1.2623&	     8.6032&	80&	     6.8\\
4&	1&	2.3226&	     6.2698&	45&	     2.7\\
4&	2&	4.8371&	     12.5925&	88&	     2.6\\
4&	3&	6.6579&	     8.2329&	45&	     1.2\\ \hline

8&	0&	8.2740&	     10.2856&	60&	     1.2\\
8&	1&	27.2746&	10.7424&	47&	     0.4\\
8&	2&	58.5001&	13.9889&	50&	     0.2\\
8&	3&	76.5998&	14.1179&	46&	     0.2\\ \hline

\end{tabular}
\caption{Performance for matrix $atmosmodd$}
\label{table-3}
\end{table}

In this experiment, the maximal speedup is 6.8 when $k$ equals zero in section three. Most speedup values are greater than one and some favorable speedups are obtained
on a GPU. Similar to the previous two experiments, both the $k$ level and the block size affect the algorithm performance on a GPU. The speedup goes down to 2.0 with
the level going up to three in section one. At the same time, the iteration number decreases from 116 to 78 which shows a better convergence improvement for a higher $k$ level.
However, the GPU time is affected and increases from 1.5976 to 3.5308, because of the additional calculations and negative influence on parallelism introduced by
fill-in positions. From the aspect of a block size, the speedup reduces quickly as a whole with a higher block size. In section four, some speedups are less than one
which means the parallel performance is exhausted by an overly large block size. Thus the GPU speed is even slower than the CPU speed under this situation.

\section{Conclusions}
\label{sec-conclusion}
We have performed a fully analysis for the decoupled block-wise ILU(k) algorithm. Its implementation and performance are investigated in detail. A few conclusions are stated below.

\begin{itemize}\itemsep1pt \parskip0pt \parsep0pt
  \item The ILU(k) algorithm can be decoupled into a symbolic phase and a factorization phase. The symbolic phase focuses on the creation of a nonzero pattern for a variable level $k$. The factorization phase assembles the lower and upper triangular matrices of ILU factorization based on the nonzero pattern created by the symbolic phase.

  \item The numerical results verify that the implementation mechanism works well for various levels $k$ and different block sizes on GPU platform. A linear solver with a block-wise ILU(k) preconditioner shows higher speed on GPU than on CPU in most situations. A favorable speedup can be acquired especially when the level $k$ and block size are lower. The speedup decreases when the level $k$ or block size increases.

  \item Although a higher $k$ is used for better convergence than a lower $k$ and a block-wise ILU(k) is more stable than a point-wise ILU(k), extra calculations decrease the parallel performance on GPU. This contributes a bottleneck for acceleration. More efforts need to be made in future study to solve this issue.
\end{itemize}
\bigskip\noindent

\textbf{Acknowledgement:}
The support of Department of Chemical and Petroleum Engineering, University of Calgary and Reservoir Simulation Group is gratefully acknowledged. The research is partly supported by NSERC/AIEE/Foundation CMG and AITF Chairs.

%-------------------------------------------------------------------------
\bibliographystyle{latex8}

\end{document}